\documentclass[twocolumn,draft,a4paper]{svjour3}[1 September 2020]
\usepackage{tikz}
\usetikzlibrary{arrows,calc,shapes,decorations.markings,scopes,patterns,trees,backgrounds,snakes}
\usepackage{oldstyle}
\usepackage[lighttt]{lmodern}
\usepackage[english]{babel}
\usepackage[T1]{fontenc}
\usepackage{graphicx}
\usepackage{epstopdf}
\usepackage{gnuplottex}
\usepackage{listings}
\usepackage{xcolor}
\usepackage{pdfpages}
\usepackage{booktabs,array}
\usepackage{amsmath}
\usepackage{amssymb}
\usepackage{microtype}
\usepackage{doi}
\usepackage{url}

\usepackage{subcaption}
\makeatletter
\newcommand\footnoteref[1]{\protected@xdef\@thefnmark{\ref{#1}}\@footnotemark}
\makeatother

\usepackage{hyperref}

\hypersetup%
{
pdftitle = {paper_wesca},
pdfkeywords = {numerics, software engineering, paper_wesca, },
pdfproducer = {xelatex},
breaklinks=true
}

\newcommand{\prog}[1]{\texttt{#1}}
\def\epem{\ifmmode{e^+e^-}\else{$e^+e^-\;$}\fi}

\title{Revisiting ``{What Every Computer Scientist Should Know About Floating-point Arithmetic}''}

\author{Vincent \textsc{Lafage}}

\institute{Vincent \textsc{Lafage}\at{}Université Paris-Saclay, CNRS/IN2P3, IJCLab, 91405, \textsc{Orsay}, \textsc{France}.\\\email{\href{mailto:vincent.lafage@in2p3.fr}{\texttt{vincent.lafage@in2p3.fr}}}}

\date{Received: \textos{2020}
}

\journalname{Computing and Software for Big Science}
\begin{document}

\maketitle

\begin{abstract}
The differences between the sets in which ideal arithmetics takes place and the sets of floating point numbers are outlined.
A set of classical problems in correct numerical evaluation is presented, to increase the awareness of newcomers to the field.
A self-defense, prophylactic approach to numerical computation is advocated.
\end{abstract}

\section{Introduction}
Precision, and its commoditization, particularly due to the plummeting of its cost and because of its reproducibility, has been the base for our industrial societies since the mid XVIIIth century. The numerical precision went even faster and farther than in mechanical manufacturing.

Thirty years after the paper
by~\textsc{Goldberg}\,\cite{Goldberg:1991:CSK:103162.103163}, the
global population of developers has arguably grown by a factor of
about fifty, and if floating-point was already deemed ubiquitous, the
then young IEEE~754 norm is now also ubiquitous, on every desk and in
every portable device. Yet ``floating-point arithmetic is still considered an
esoteric subject by'' numerous computer practitioners.

While numerical computations are pervasive nowadays, they may appear for what they are not. Obviously, one needs computers, languages, compilers, and libraries for computing. Furthermore one has to rely on some standard, such as IEEE~754~\cite{Kahan:1979:PFS:1057520.1057522,30711,27840} 
and its evolutions~\cite{4610935,8766229}, that quickly became ISO standards~\cite{IEC559:1989,ISO60559:2011,ISO60559}. We note that this standard emerged in a rich environment of incompatible formats, and is now supported on all current major microprocessor architectures. Yet all this is but modeling the numbers. Our familiarity with numbers typically used for computational purpose, \emph{i.e.\ }the real numbers, may then blur our understanding of floating point numbers.

We will first discuss the common floating point representation of real numbers in a computer and list the differences between the properties of the numbers and the properties of their representation. We will then illustrate common examples where a naive implementation of the computation can lead to unexpectedly large uncertainties. In the next part, we delineate some approaches to understand, contain or mitigate these accuracy problems. We then conclude on the need for a broader awareness about floating point accuracy issue.

\section{From real numbers to floating point numbers}
Our ``intuitive'' understanding of numbers, for measurements, payments, arithmetics and higher computations, was developed over years of schooling and essentially builds on decimals which is not only a base-ten positional numeral system,
but also the set of all decimal fractions, called $\mathbb {D} =\left\{{\frac {n}{10^{p}}}, n\in \mathbb {Z}, p\in \mathbb {N}\right\}$ for decimal.
This decimal notation, obvious for us, is
rather new, dating back to François~\textsc{Viète} in the XVIth century, who introduced systematic use of decimal fractions and modern algebraic notations,
and John~\textsc{Napier} in the early XVIIth century (improving notations of Simon \textsc{Stevin}) who popularized modern decimal separators (full stop dot as well as comma) to separate the integral part of a number from its fractional part, in his \emph{Rabdologia} (\textos{1617}).
Decimal notation being central for the metric system, its use is widely spread.

But computing with machines rely essentially on binary representations, and computers are using subsets of the binary fractions, called $\mathbb {B} =\left\{{\frac {n}{2^{p}}}, n\in \mathbb {Z}, p\in \mathbb {N}\right\}$ for binary.

$\mathbb{B}=\mathbb{Z}[1/2]$ is generated by adjoining $1/2$ to $\mathbb{Z}$, in the same way that $\mathbb{D}=\mathbb{Z}[1/10]$.
This construction is a rational extension of ring $\mathbb{Z}$ or localization of the ring of integers with respect to the powers of two.
Therefore, algebraically, $(\mathbb{B}, +, *)$ is a ring (like $\mathbb{D}$).
As it is not a field, division, as well as all rational functions, will require rounding.
As $\mathbb{B}$ (like $\mathbb{D}$) is only a small subset of the algebraic numbers, all roots will require rounding.
As $\mathbb{B}$ (like $\mathbb{D}$) is only a small subset of the real numbers ($\mathbb{R}$), all transcendental functions\footnote{\label{note_transcend}as opposed to \emph{algebraic} functions, \emph{i.e.\ }functions that can be defined as the root of a polynomial with rational coefficients, in particular those involving only algebraic operations: addition, subtraction, multiplication, and division, as well as fractional or rational exponents} will also require rounding (even elementary ones: trigonometric functions, logarithms and exponential ; and special functions such as \textsc{Bessel} functions, elliptic functions, \textsc{Riemann} $\zeta$ function, hypergeometric series, polylogarithms\dots).
As $\mathbb{B}$ is a ring, we get most usual good properties we are used to in mathematics:
\begin{itemize}
\item commutativity:\\
  $\forall{}x\in\mathbb{B},\forall{}y\in\mathbb{B}\quad x+y = y+x$,\\
  $\forall{}x\in\mathbb{B},\forall{}y\in\mathbb{B}\quad x\times{}y = y\times{}x $
\item associativity:\\
  $\forall{}x\in\mathbb{B},\forall{}y\in\mathbb{B},\forall{}z\in\mathbb{B}\quad x+(y+z) = (x+y)+z$,
  $\forall{}x\in\mathbb{B},\forall{}y\in\mathbb{B},\forall{}z\in\mathbb{B}\quad x\times{}(y\times{}z) = (x\times{}y)\times{}z$
\item distributivity:\\
  $\forall{}x\in\mathbb{B},\forall{}y\in\mathbb{B},\forall{}z\in\mathbb{B}\quad x\times{}(y+z) = x\times{}y+x\times{}z$
\end{itemize}
Unfortunately, in practice we are limited to elements of $\mathbb{B}$ with a finite representation.

We call these ``binary\footnote{We will not deal here with other bases, except to illustrate the cancellation phenomenon and double rounding in the more natural base ten} floating point numbers'', with a given number of \emph{digits}, in this case \emph{bits}: $m$ bits for their \emph{mantissa} or \emph{significand}, and $e$ bits for their \emph{exponent}. Let's note these sets $\mathbb{F}_{m,e}$: they build an infinite lattice of subsets of $\mathbb{B}$. The size $e$ of the \emph{exponent} and $m$ of the \emph{mantissa} includes their respective sign bits. Furthermore some implementation do not store the most significant bit of the mantissa as it is one by default, then saving one bit for accuracy. Note that, this saving is possible only as we use a base two representation. Most of the time we rely on the hardware embedded $\mathbb{F}_{24,8}$ for single precision and $\mathbb{F}_{53,11}$ for double precision.

Topologically, while $\mathbb{B}$ (like $\mathbb{D}$) is dense in $\mathbb{R}$, which would allow for arbitrarily close approximations of real numbers, the finite counterparts of $\mathbb{B}$ (or $\mathbb{D}$) used for floating points, are nowhere dense sets.

Similarly, while $\mathbb{B}$ (like $\mathbb{D}$) are closed under addition and multiplication, their finite counterparts $\mathbb{F}_{m,e}$ is not.
\begin{itemize}
\item for multiplication, the number of significant digits of the product is usually the sum of significant digits of factors.
\item for addition, the sum will be computed after shifting the less significant term relative to the most significant one, the exact resulting significand getting larger than the original one.
\end{itemize}
If these results happen to fit into the fixed size, then the computation has been performed exactly, else an \texttt{inexact} flag can be raised by the processor and this requires rounding.
This is a source of trouble: we cannot catch, or \emph{trap} in the specific case of IEEE~754, all \texttt{inexact} exception that would appear in any realistic computation. Most of them even appear in the simplest operation (addition), and then we do not have a choice to handle them: dynamically adaptable precision is possible only with a massive loss of speed, and is useless as soon as roots or transcendental functions\footnoteref{note_transcend} are used. These exceptions will appear for all level of precision.

When representing a real number $x$, one find the two nearest representable floating point numbers: $\mathtt{x_-}$ immediately below $x$ and $\mathtt{x_+}$ immediately above $x$.
Let us note $m$ the middle of the interval $[\mathtt{x_-},\mathtt{x_+}]$ (a real number that cannot be represented in this floating point format).
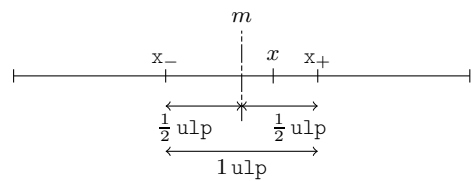
\begin{figure}[h!]
  \begin{center}
    \begin{tikzpicture}[scale=2]
      \draw [dash pattern={on 7pt off 1pt on 2pt off 1pt}] (0.5,0.2) -- (0.5,0.8) node[above]{$m$};
      \draw (0.707,0.45) -- (0.707,0.55) node[above]{${x}$};
      \draw (0,0.5) node[above]{$\mathtt{x_-}$};
      \draw (1,0.5) node[above]{$\mathtt{x_+}$};
      \draw[<->] (0,-0) -- (1,-0) node[pos=0.5,below]{$1\,\texttt{ulp}$};%
      \draw[<->] (0,0.3) -- (0.5, 0.3) node[pos=0.25,below]{${1\over2}\,\texttt{ulp}$};%
      \draw[<->] (0.5,0.3) -- (1, 0.3) node[pos=0.75,below]{${1\over2}\,\texttt{ulp}$};%
      \draw (-1,0.5) -- (2,0.5);
      \draw[snake=ticks,segment length=2cm] (-1,0.5) -- (2.01,0.5);
    \end{tikzpicture}
  \end{center}
  \caption{The basics of rounding}
  \label{fig:rounding}
\end{figure}
In order to measure the effect of rounding, one uses the smallest discrepancy between two consecutive floating point numbers, $\mathtt{x_-}$ and $\mathtt{x_+}$: they differ by one unit in the last place (\texttt{ulp}): see Fig.~\ref{fig:rounding}. So when rounding to the nearest, the difference between the real number and the rounded one should be at most half a \texttt{ulp}: this is called \emph{correct rounding}. As this proves difficult, a less strict requirement is to guarantee the rounding within one \texttt{ulp} (\emph{i.e.\ }all bits correct but possibly the last one), which is called \emph{faithful rounding}.
In Fig.~\ref{fig:rounding}, both $\mathtt{x_-}$ and $\mathtt{x_+}$ are \emph{faithful roundings} of $x$, but only $\mathtt{x_+}$ is the \emph{correct rounding} of $x$.

Thus, the behavior of floats is not the same as the behavior of binary numbers.
Addition and multiplication inside the set of binary numbers are stable while for floats,
it will be inexact most of the time, requiring round-off.
From an object oriented point of view, while the floats $\mathbb{F}_{m,e}$ being smaller sets should be subtypes of the binary number as a subset of $\mathbb{B}$,
the extra amount of exceptions they bring come in violation of the \textsc{Liskov} substitution principle~(LSP)~\cite{10.1145/62138.62141,10.1007/3-540-47910-4_8,10.1145/197320.197383}\footnote{\textsc{Liskov}'s notion of a behavioral subtype defines a notion of substitutability for objects; that is, if S is a subtype of T, then objects of type T in a program may be replaced with objects of type S without altering any of the desirable properties of that program (e.g.\ correctness).}.

We must not regard these subsets as behavioral subtypes, which goes against our intuition.

In the same way, while the defined exceptions for single and double are the same,
they will not be triggered as often for the highest precision,
and in particular for higher precision representation of the same pure number:
a product of singles can overflow, where the product of the same values would be computed correctly,
and even exactly, with doubles.
So we do not get a tower or lattice of inheriting types, the smaller precision deriving from the higher ones,
as new behavior (in the form of more exceptions) would emerge from the subtype, in violation of the LSP{}:
in particular, precondition on domain of functions are stricter in the subtype (to prevent overflow or underflow).

Some may consider that subtypes extends the base type, with extra information.
In this respect, higher precision would be considered as subtypes of lower precision,
extending it with extra digits, and keeping the same methods.
Yet, the exception being triggered more often with the lower precision,
the higher precision subtype would lead to different behavior (at least better correctness) of the same program.

In order to avoid too many exceptions, particularly exception out of closure, extra elements are added to the set.
First, to take overflow into account, $\mathbb{B}$ is affinely extended by signed infinities: $\overline{\mathbb{B}} = \mathbb{B}\cup\{-\infty, +\infty\}$.
Then signed zeros $(-0,+0)$ are included for consistency. One more element is convenient to cover for undefined results such as functions on real numbers returning complex numbers, but also indeterminate forms: \texttt{NaN} for ``Not a Number''.
With this we turn our number system into a broader representation. We then get another source of \texttt{NaN}, as the result of operations involving \texttt{NaN}.
Which in turns represents different cases so that we may want different \texttt{NaN} in our representation to give \texttt{NaN} a more specific meaning: these extra informations can be stored in the \texttt{NaN} sign bit and in a payload of $23$ bits for single precision (or $52$ bits for double)\footnote{as is the not well known and even less used case in IEEE~754. Language constructs or libraries to capitalize on this are unfortunately absent. As these payloads are not well known, standard libraries avoid the tedious work of correctly propagating them, which limits even more their usefulness.}.

With these extensions the lack of closure of bare float numbers is solved. But this solution comes at the cost of losing associativity and distributivity (a major behavioral shift), that are hard-wired in our ways of thinking about numbers, in the same way the implicit hierarchy of precisions is. Beyond algebraic properties, \texttt{NaN} also destroys the total order that exists in $\mathbb{R}$ and $\mathbb{B}$. For low-level programming and optimization, it also breaks the connection between logical equality and binary representation equality.

While the exceptional behaviors are usually apparent, even when replaced by their representations \texttt{NaN} or \texttt{Infinity}, and the \texttt{inexact} exception can be triggered with only half a \texttt{ulp} of loss of precision, a much more problematic case will not trigger any more specific exception:
``\emph{catastrophic cancellation}''~\cite{Goldberg:1991:CSK:103162.103163}, also known as ``\emph{loss of significance}''
occurs when an operation on two numbers increases relative error more than it increases absolute error,
for instance, when subtracting close numbers. As a result of this closeness, the number of significant digits in the outcome is less than expected with the precision.
For instance, using here decimal floating points with three digits mantissa for the sake of clarity, when one considers $a=3.34$ and $b=3.33$, on the one hand the mathematical difference $a-b=0.01$ is computed exactly with the given floating point format. This computation, while exact, is a cancellation: one recovers only one significant digit out of three. As this result could not be more exact, this cancellation is said to be \emph{benign}. On the other hand, the mathematical difference of the squares $a^2-b^2=0.0667$ is $6.67\times10^{-2}$ in the given floating point format, but will be evaluated naively as $0.1=10^{-1}$ with the rounding. No digit of the floating point result is correct, and this result differ by $50\%$ relative error or 333 units in the last place: this cancellation is \emph{catastrophic}.
There is no systematic way to avoid catastrophic cancellation, or to know how many digits have been lost. In particular, catastrophic cancellation will not always lead to a zero, as illustrated with this example.

A corner case among the already pathological behavior of catastrophic cancellation occurs when the absolute value of the computation is below the smallest normal floating point number, \emph{i.e.\ }$2$ to the power of the smallest exponent for the given format: $2^{-128}\sim 0.29\times10^{-38}$ for single precision, $2^{-1024}\sim 0.56\times10^{-308}$ for double precision.
A strict application of the representation would lead to an underflow exception, turning the result to $0$, even if the mathematical result is not strictly null. This could occurs when subtracting close small numbers. To prevent this total loss of relative accuracy, a special representation is included for \emph{subnormal} numbers (formerly known as \emph{denormal} numbers). While these numbers provide algorithms with more stability, dealing with this anomalous representation can slow down computation by as much as $100$: this fact brands this part of the norm a shame to high performance computing.
A first thing to take into account is that such a slowdown points out a likely improper normalization of the values (see subsection~\ref{sub:nondimensionalization}). A second aspect when subnormal numbers occur, is that extending the precision from single to double is likely to cure the problem at a vastly less cost than redesigning the algorithm or renormalizing the values.

\section{Use cases}
In this section we will present a collection of use cases where the developers should be aware of the inaccuracies of the floating point computation, and refrain from using the naive approach.

\subsection{Converting decimal to binary and binary to decimal}
As $\mathbb{B}\subset\mathbb{D}$, every binary floating point number has an exact decimal representation within a finite number of decimal digits. But not conversely: the decimal numbers with which we count our cents may not always get a finite binary representation to the disarray of beginners, the archetypal example being the violation of the mathematical identity $0.1 + 0.2 = 0.3$, {whatever the floating-point precision, but more evidently in single precision}.

So when, frequently, there is no $n$ digits binary number corresponding exactly to the input decimal number, we round off towards the closest $n$ digits binary number. This requires converting decimal to a precision larger than $n$ binary digits to correctly round it. Hidden in most hardware, below the least significant bit, one finds first a guard bit, a round bit and then a sticky bit: the two first are just two extra bits of precision, while the sticky bit keeps track of all non zero bits lost during the right-shifts of the mantissa.
But the decimal to binary conversion is affected by a problem related to the table-maker dilemma~(TMD)~\cite{castiel:inria-00000567,10.1117/12.505591} and these three extra bits may not be sufficient to proceed to the correct rounding.
The system standard \texttt{C} libraries usually provide us with elementary and approximate versions of functions~\texttt{atof} or \texttt{strtof}.

As for the binary to decimal conversion, although every binary floating point number has an exact and finite decimal representation, it is usually not this exact conversion that is used in practice, for instance in a \texttt{printf}. This faithful translation would require more digits than the effective precision used for the floating-point format, leading to an absurd overprecision.
The decimal number that should be returned is then not the exact one, but the decimal with the shortest expression that stand within one half \texttt{ulp} of the binary number. This minimization leads to comparisons and computations with precision higher than the target precision for the routine. The system libraries typically provided a more pragmatic approach, not aiming at correctly rounded but \emph{good enough} results~\cite{QuickDirtyToDecimal}.

Yet, correctly rounded conversion routines have existed
for the last thirty years~\cite{Gay90correctlyrounded}~({\texttt{C}~source available on~\href{http://www.netlib.org/fp}{\texttt{www.netlib.org/fp}}}). Faster algorithms even appeared~\cite{Loitsch:2010:PFN:1809028.1806623}. But they did not seem to make it in our compilers and system or language libraries. While the compilers have gone a long way to evaluate the literal and named constants in extra precision, and that mitigates the problem, the bulk of the conversions comes from decimal data input/output: as long as these are not commoditized as language construct, these aspects of round-off are often considered for experts.

This kind of round-off can pile-up in loop of conversions (from one program to another, for instance when storing intermediate results in decimal format). As a consequence, one should avoid decimal storage as much as possible and prefer binary storage, not only for performance, but also for precision.
Moreover, the number of bits for accurately encoding a given number of decimal digits is not the naive estimation: it has been known for over fifty years~\cite{10.1145/363067.363112}, that the $24$ bits of single precision are not enough for an accurate representation of $7$ decimal digits, while $24\log_{10} (2)>7$. Hence caution is needed when estimating the required precision, usually on number of decimal grounds.

\subsection{Double rounding}\label{sub:double_rounding}
The classical issue of double rounding occurs when using two different precisions: a high precision and a lower one. The exact result is first rounded to the high precision, then to the lower one. For instance with base 10 floating points, we will consider a mantissa of 2 digits for the high precision and 1 digit for the lower one: with an exact result of $3.47$, the rounding to the high precision gives $3.5$, and the rounding to the lower precision gives $4$, whereas the direct rounding of the exact result to the lower precision is $2$. 
Double rounding happens when computing in extended precision (such as the hardware implementation of $\mathbb{F}_{64,16}$) that one converts to the smaller working precision: even if the extended precision has been correctly rounded (within half a \texttt{ulp}), the rounding to working precision may lead only to a faithful rounding (within one \texttt{ulp}). The correct solution is addressed in Ref.~\cite{Normale04whendouble}. In particular, while double correct rounding is not always correct, the double faithful rounding is  faithful.

\subsection{Difference of two squares}\label{sub:factorization}
The difference of two squares $a^2-b^2$ seems innocuous, and if the units of the program have been chosen correctly, we can easily avoid the overflow risk (however it is enough to trigger this overflow that $a$ or $b$ be of order of $10^{19}$ in single precision instead of $10^{154}$ in double precision). The real issue is, when the two terms are close, in the possibility to exhibit ``\emph{catastrophic cancellation}''~\cite{Goldberg:1991:CSK:103162.103163}, also known as ``\emph{loss of significance}''.
These cancellation occurs in an operation (usually a subtraction) on previously truncated operands, for instance, the product of other operands (here, squares).
The numerical result then has only a limited precision compared to that with which the operand computations were made. Even the sign of the result can be wrong.

\textit{A contrario} the same difference on exact operands leads to a benign cancellation.

Using double precision is usually the way to mitigate this effect.

The best solution is to use a factorized formulation of the computation $a^2-b^2 = (a-b)\cdot(a+b)$. The remaining subtraction can be the cause of ``\emph{benign cancellation}'' at most.
In this case factorization used to provide a speed gain (a single multiplication instead of two) on top of the precision gain\footnote{technically, when $a$ and $b$ are of different magnitude, the non-factorized version suffers one less rounding error, therefore is more precise, but without affecting the meaning of the result}. On contemporary architectures, the speed gain has disappeared, as multiplication is now as fast as addition.
In this simple case, the factorization is trivial, but in the expression of the area of a triangle, or the volume of a tetrahedron, the work required is arduous. More generally, in the determinantal expressions common in geometry, there is no guarantee of the existence of these factorizations, and when they exist such results seems to be found by sheer luck.

\subsection{Summation}\label{sub:sum}
Another well known algorithm among the community is Monte-Carlo integration: basically, it is a loop sampling a configuration space, then accumulating a given contribution as function of the sampled point.
We therefore proceed to add fluctuating values in a growing accumulator.
The more one accumulates, the more the accuracy of the individual contribution is blurred.

Ideally, one should sort the contributions from the smallest to the largest ({for contributions with the same sign}) to be sure to lose as little precision as possible. One could also apply a ``divide and conquer'' algorithm such as pairwise summation~\cite{doi:10.1137/0914050}. Another solution is to accumulate in extended arithmetic, or to use the summation algorithm of \textsc{W.~Kahan}~\cite{Kahan:1965:PRR:363707.363723} or other compensated sums~\cite{10.1145/980175.980177,doi:10.1137/1037130,babuska1968numerical,doi:10.1137/0914050,doi:10.1137/030601818,doi:10.1137/07068816X}.

Let us specify this essential aspect of sums which is the accumulation of truncation error.
As this error will statistically be distributed uniformly as much upwards as downwards, one can assimilate it to a random walk, and in the hypothesis of terms of comparable size, to expect for the sum of $N$ terms to a resulting absolute uncertainty on the total of the order of $\mathcal{O}\,\left({\sqrt{N}}\right)$.
This does not seem alarming since the associated relative uncertainty will be in $\mathcal{O}\,\left({1/\sqrt{N}}\right)$.
This approach is however naive since the absolute error of truncation of a sum is like the relative error compared to the largest term \emph{i.e.\ }the accumulator.
The relative error varies like $\mathcal{O}\,\left({\epsilon\sqrt{N}}\right)$. Comparatively, pairwise summation has a much tamer relative error varying like $\mathcal{O}\,\left({\epsilon\sqrt{\log_2{N}}}\right)$.
As a matter of fact, the naive addition of $4\times10^6$ contributions in a half precision accumulator causes a relative uncertainty of $100\,\%$, and in the case of a single precision accumulator, it reduces the precision within four significant digits at best (if each contribution is correctly calculated to the accuracy of the machine and correctly rounded). The average Monte Carlo can easily produce such a four million events sampling.

We can also combine these mild (and all the more treacherous) losses of precision with catastrophic cancellation. This is, for instance, the case of series falling under the \textsc{Leibniz}'s test for alternating series: a sum of terms of alternating signs, of zero limit at infinity and of decreasing absolute values. The convergence of these series is mathematically guaranteed, and the uncertainty on the partial sum is bounded from above by the absolute value of the first neglected term in truncating the serie. Each pair of successive terms of these series is thus concerned with a potentially disastrous compensation. The latest can be overcome in the case where an analytical expression of terms exists, but there is no general solution for purely numerical expressions.

\subsection{Quadratic formula}\label{sub:quadratic}
Let's start with the famous quadratic equation
\begin{eqnarray}
  {\displaystyle ax^{2}+bx+c=0,\quad{\text{where}}\quad{}a\not=0,}\\
  {\displaystyle{\text{with}}\quad\Delta={b^{2}-4ac},}
\end{eqnarray}
its corresponding discriminant, we note, for $\Delta>0$, its two exact solutions:
\begin{eqnarray}
{\displaystyle x_\pm={\frac {-b\pm {\sqrt\Delta}}{2a}}.}
\end{eqnarray}
We can identify two possibilities for ``\emph{catastrophic cancellation}'': the first, in expressing the discriminant, and the second in the compensation between $-b$ and $\sqrt\Delta$. The latter can be overcome reformulating the roots as:
\def\sgn{\mathop{\mathgroup\symoperators{}sgn}\nolimits}\\ 
\begin{eqnarray}
q   & = & {-b-\sgn(b) {\sqrt\Delta}} = {-\sgn(b)\left({|b| + {\sqrt\Delta}}\right)}\\ 
x_1 & = & {\frac {q}{2a}},\\
x_2 & = & {\frac {2c}{q}} = {\frac {c}{ax_1}}.
\end{eqnarray}
$q$ is then always a sum of two like-sign terms, without catastrophic cancellation.
We thus obtain an accurate expression without execution branching.
On the other hand, avoiding execution branching blurs the link between the two pairs of roots, and a test might be necessary to match $x_1$ to $x_+$ or $x_-$.

A solution for taming the catastrophic cancellation in the discriminant will be presented in Subsection~\ref{subsub:fma}.

\subsection{Polynomial evaluation}
It is well known that polynomial function should not be evaluated in the naive way: on the one hand, because of the computational cost of all the exponentiation, on the other hand because of the accuracy loss it represents. As we have seen in Subsection~\ref{sub:factorization}, when a factored out form exists, it had better be exploited. In addition, when factorization is not known, the best approach is \textsc{Horner-Ruffini} method~\cite{10.2307/107508}:
\begin{eqnarray}
{\begin{aligned}p(x)&
=a_{0}+a_{1}x+\cdots+a_{n-1}x^{n-1}+a_{n}x^{n}\\&
=a_{0}+x{\Big({a_{1}+x{\big({\cdots +x(a_{n-1}+x\,a_{n})\cdots}\big)}}\Big)}\\\end{aligned}}
\end{eqnarray}
It is not only a gain in speed, through the saving of operations that it achieves, but also in accuracy, partly for the same reason, and more importantly it is a guarantee of stability of the result and safety against intermediate overshoots, particularly in the case of simple precision. It is even possible to increase the evaluation speed~\cite{Knuth:1962:EPC:355580.369074}.

The democratization of ``multiply-accumulate'' machine instructions (\texttt{fma}) combining addition and multiplication (see Subsection~\ref{subsub:fma}), leads us to strongly seek their use in this case~\cite{Graillat:2006:ICH:1141277.1141585}.

Furthermore, a polynomial being a sum (see Subsection~\ref{sub:sum}), compensation summation techniques, such as the summation algorithm of \textsc{W.~Kahan}~\cite{Kahan:1965:PRR:363707.363723}, must be taken into account.

\subsection{Variance estimate}
The computation of variance estimator is a canonical example of ``\emph{catastrophic cancellation}'' in that it is a difference of squares. Let us consider the variance of a sample of $N$ individuals.
\begin{eqnarray}
\sigma^2={\overline{{(x-{\bar {x}})}^2}}={\overline {(x^2)}}-{\bar {x}}^2= {\frac {\sum _{i=1}^{N}x_{i}^2-\frac {{(\sum _{i=1}^{N}x_{i})}^2}{N}}{N-1}}.\!
\end{eqnarray}
A golden rule for students used to be ``\emph{never compute the variance with single precision}''\dots{}
Because the variance so computed could become negative, it appears abruptly when one extract the square root to obtain the standard deviation. As evaluation occurs more and more online, consequences can be dire if proper safe-guard exception handling is not provided.

Certainly, a two-pass approach, when possible, allows to stick with single precision.
But one can also take advantage of~\textsc{Welford}~\cite{doi:10.1080/00401706.1962.10490022} one pass online algorithm. The cost for accuracy then comes as an extra division for each iteration loop.

The problem also appears when evaluating covariance matrices and in evaluating centered moments of higher orders.

\subsection{Area of a triangle}
The formula of \textsc{Heron of Alexandria}, allows to calculate the area $S$ of any triangle knowing only the lengths $a$, $b$ and $c$ of its three edges:
\begin{eqnarray}
 \displaystyle S={\sqrt{p(p-a)(p-b)(p-c)}}\quad
\end{eqnarray}
 with\quad $p={\frac{a+b+c}{2}}$\quad{the half perimeter}.
This nice symmetrical formula exhibits an instability during the numerical computation, which appears for a needle-like triangle, \emph{i.e.\ } a triangle of which one edge is of particularly small dimension compared with the others (confrontation of small and large values).

By labelling the names of edges so that $a>b>c$, and by reorganizing the terms to optimize the quantities added or subtracted, W.~\textsc{Kahan} proposed a more stable formula for $S$ in Ref.~\cite{Kahan_Triangle}:
\begin{eqnarray}
{\displaystyle {\frac{1}{4}}{\sqrt{\left[{a+(b+c)}\right]\,\left[{c-(a-b)}\right]\,\left[{c+(a-b)}\right]\,\left[{a+(b-c)}\right]}}}\cr
\end{eqnarray}
While its apparent symmetry is lost, this formula is much more robust than the naive one.
He similarly provided a stable formula for the computation of the volume of a tetrahedron based on a non trivial factorization~\cite{Kahan_Tetrahedron}.

These calculations may seem hazy textbook case: they are in fact only the expression with small dimensions ($n=2$, $n=3$) of factors which intervene in more general geometric calculations, such as the relativistic phase space volume~\cite{byckling1973particle}. These factors turn out to be determinants, notoriously delicate to evaluate precisely if one proceeds naively\footnote{Multivariate polynomial in nature, while relying on ``simple'' arithmetics, these multi-linear alternating form are the best place to exhibit massive cancellation. Computation examples defeating even quadruple precision rely on these higher degree polynomial cancellation (see Ref.~\cite{RUMP1988109,loh2002rump})}. For instance, if we study the production of intermediate bosons and electrons, their square masses are in a ratio of $10^{10}$ and certain kinematic factors will be analogous to the surface of a triangle with edges in such a ratio.

\subsection{Complex Arithmetic and Analysis}
Complex numbers are a powerful tool for the physicist and the engineer. Unfortunately, they are often presented as the simplest example of Cartesian coordinates, which often gives rise to high school-like implementations\footnote{see comments such as \texttt{// XXX:\@ This is a grammar school implementation.} in the \texttt{complex} header of \prog{g++}} which are mathematically correct but numerically fragile. These expressions are notably liable to involve intermediate overflow although their result is representable (complex modulus, complex division and even complex multiplication). For instance, when naively evaluating the modulus of a single precision complex number with the largest component greater than $2\times10^{19}$, the intermediate sum of squares will overflow while the resulting modulus is representable as a single precision floating point number.
Other languages integrate complex numbers better, such as \texttt{Fortran}, \texttt{Ada} and \texttt{D}, and even better, \texttt{C99}.

As for complex analysis, particularly when evaluating radiative corrections in particle physics, but also when modeling 2D irrotational fluid flows, the subtleties of managing the zero sign in the vicinity of the cuts associated with connection points, can lead to results which are not trivially false~\cite{Kahan:CSD-92-667}. It is therefore important to use established libraries to benefit from guaranteed values, for example for one-loop functions in radiative corrections~\cite{vanOldenborgh:1990yc}.

\subsection{Function evaluation}
Function evaluation often start with a range reduction: trigonometric functions are first linked to their expression in the range $[0;{\pi\over4}]$, and building upon the binary representation of the argument, square root over $[{\sqrt2\over2};\sqrt2]$, exponential in the range $[0;\ln2]$, logarithm over $[-\sqrt2;\sqrt2]$\dots
This apparently trivial transformation can in itself require extra computations to be really accurate.
But before algorithms such as power series expansion or \textsc{Newton}-like iteration can be used efficiently, an extra step needs to be taken, as these reduced intervals are already too large for quick convergence. This is commonly achieved through lookup tables (LUT).

\subsection{Function Minimization}\label{sub:minimization}
 This almost daily activity in our community, seems well under control. If we look more closely, in the vicinity of its minimum, the studied function is the sum of the minimum value and of a quadratic form (\textsc{Taylor}'s theorem). With this form's coefficients close to unity, a relative displacement of $\mathcal {O}\,(h)$ around the position $x_0$ of the minimum will generate a relative displacement of $\mathcal{O}\,(h^2)$ of the function to minimize. However, the smallest relative measurable displacement being bounded by the \emph{machine epsilon} $\epsilon=2^{-m}$, the variation $\delta{}x_0$ which will not generate a perceptible difference will be $\mathcal{O}\,\left({\sqrt{\epsilon}}\right)$, half the machine precision. Therefore if we calculate the function with double precision, the uncertainty on the position of the minimum corresponds to the single precision. If we now evaluate the function to minimize in single precision, the position of the minimum will only have the resolution of the half precision. Moreover if we estimate the function by means of the half precision in vogue, the position will be known with less than two significant figures. This half precision being supposed to serve the training in \emph{machine learning}, which rely essentially on minimization, using half precision casts some doubts on the accuracy of the results. 

Some studies have tested the consequences of using half precision in training neural networks, including  the use of mixed precision with single precision accumulators (see Refs.~\cite{courbariaux2015,DBLP:journals/corr/abs-1710-03740} or fixed point number instead of floating point numbers (see Refs.~\cite{gupta2015deep,DBLP:journals/corr/abs-1809-00095}). These precisions are sufficient for training and running the network. The question of the accuracy of a neural network reproducing a numerical function is yet open. 

\section{Some approaches to control floating point uncertainty}
Now cautious that our prejudice will work against us, we need some way to compute nevertheless.

\subsection{Condition number and stability}
As we have seen, some computations have a risk of developping large relative uncertainties even with precise input. In this subsection, we will attempt to quantify the amplification factor between input uncertainty and output uncertainty, and to discriminate the part inherent to the function under study, called the condition number, from the part caused by the evaluation algorithm.

The condition number $\kappa$ inherent to a mathematical function $f$ evaluated at argument $x$ with approximate value $x_a$ compares the relative uncertainty on $f (x)$ to the relative uncertainty on the argument. It takes the form:
\begin{eqnarray}
\kappa\,(x)
=\frac{\left|{\frac{f(x_a)-f(x)}{f(x)}}\right|}{\left|{\frac{x_a-x}{x}}\right|}
=\frac{\left|{\frac{f(x_a)-f(x)}{(x_a-x)}}\right|}{\left|{\frac{f(x)}{x}}\right|}
\sim\left|{\frac{x{f'(x)}}{f(x)}}\right|
\end{eqnarray}
$\kappa$ is also known as the absolute value of the \emph{elasticity} of the function $f$.
Note that $\kappa$ is dimensionless.
This simple expression reveal that even the most elementary function is prone to large uncertainties: considering the addition of a constant $c$ ($f:x\to{}x+c$), the corresponding condition number $\kappa=\frac{x}{x+c}$ exhibits a singular behavior around $x=-c$, while being tame on the vast majority of value: next to the cancellation point, relative error on the result will be far larger than relative error on the argument. This particular example illustrates that catastrophic cancellation is amplified by a large condition number even before rounding takes place.

Beyond this formal definition, a more concrete consequence of the condition number will help us measure its value as a tool: as the smallest relative uncertainty achievable with a given precision is of half a \texttt{ulp}, a neatly computed parameter of a function $f$ will produce a relative uncertainty of $\frac{\kappa}{2}$\,\texttt{ulp} in the result of $f$. Hence $\log_{10}\,\kappa$ is the number of digits of accuracy lost in this computation, and $\log_2\,\kappa$ is the number of bits of accuracy lost in this computation.

Moreover, the composition of functions will lead to a net condition number equal to the product of component functions condition numbers. The $23$ bits of single precision can then succumb to repeatedly composed evaluation of a function with as moderate a condition number as $2$.

This points out the need for different levels of cautiousness in the evaluation strategy of a function, depending on the argument: some subdomains of the definition domain do not exhibit large condition number, others need more care. {From the definition, one can easily demonstrate that only power laws $x\to{}C\times{}x^n$ (with $C$ and $n$ some real constants) have a constant condition number $\forall{}x,\,\kappa\,(x)=n$}. Hence even functions with a nice unique expression may need piecewise evaluation. Beside lacking of aesthetic appeal, this casuistry has a clear impact on performance as the required tests potentially lead to branch misses. Moreover, these tests come on top on those needed for checking that the argument is inside the function's domain of definition.

The condition number of the problem reveals that the precision of the computation may have to be larger than the target accuracy.

Note that until this point, we have considered only exact mathematical expression: the condition number $\kappa\,(x)$ does not consider approximations and rounding. Further uncertainties come from the numerical algorithm used to evaluate the expression: approximation, for instance truncation of a power or \textsc{Fourier} series expansion, as well as the rounding, on each term and on the sum. The condition number and extra uncertainties contribute to define an overall \emph{stability} of the algorithm under small perturbations. These concepts are as old as computers~\cite{10.1093/qjmam/1.1.287} but deserve to be given a larger exposure: these were for instance not described in Ref.~\cite{Goldberg:1991:CSK:103162.103163}, despite the numerous examples of this paper being developped around singular points for their respective condition number.

Once the singular values of the condition number is found (most of the time simple zeros of function $f$ will be singular for its condition number $\kappa$), an \emph{asymptotic analysis} of $f$ around these singular values must be performed to assess the more stable expression: for instance, considering $f=x\to\ln\,x$, one finds that $\kappa\,(x)=\frac{1}{\ln\,x}$, with singular value $x=1$, which brings our attention to the expression of $f$ around $1$, $f\,(x=1+h)=\ln\,(1+h)$ for which the condition number becomes:
\begin{eqnarray}
\kappa\,(h)=\frac{h}{(1+h)\ln\,(1+h)}\underset{h\to0}{\sim}\frac{1}{(1+h)}
\end{eqnarray}
which is no more singular for $h$ near $0$ (but is singular for $h$ near $-1$). Yet as the expression involves an indetermination, it needs to be evaluated with care as described in Ref.~\cite{Goldberg:1991:CSK:103162.103163}. Nowadays, the mathematical libraries provide us with \texttt{log1p}.

In real situations, the condition number must be extended beyond univariate function by using an arbitrary norm instead of absolute value, and one must take into account the uncertainties on physical constants as possibly greater than half a \texttt{ulp} (as mathematical constants have zero for their condition number).

\subsection{Nondimensionalization}\label{sub:nondimensionalization}
Scientists often compute to solve equations or simulate phenomena where variables and parameters have specific units and physical dimensions, following dimensional analysis. Dimensional homogeneity requires that all analytic functions and their arguments must be dimensionless, as these functions can be expressed as an infinite power series of the argument and all terms must have the same dimension.
Usage is to extract the \emph{characteristic units} of the problem at hand, and use these for \emph{rescaling} the dimensional quantities to dimensionless ones (\textsc{Vaschy}-\hspace{0pt}\textsc{Buckingham}-\hspace{0pt}\textsc{Bertrand} $\pi$ theorem). This choice of units allows one to simplify the equation, with parameters both fewer and dimensionless (such as the \textsc{Reynolds} number in hydrodynamics, among others).
Beside this aspect, useful to understand the problem, it also proves beneficial numerically as the dimensionless variables are of the order of unity, where the density of representable floating points is the highest (and where the overflow and underflow are as far as possible).

\subsection{Interval arithmetic}
If one would like a higher precision than the native types of the processor, there are different solutions such as double-double precision, quadruple precision, arbitrary precision, rational representations\dots

An important alternative is to represent each variable by a pair of floating point numbers representing the bounds of an interval ensuring that it contains the true value. Numbers arithmetic then becomes interval arithmetic, which is more than twice as expensive to evaluate, but provides mathematical certainties. This approach can be traced back to \textsc{Archimedes} and appeared quickly but discreetly in the history of computers. Beyond the calculation of each terminal, it is advisable to round the lower limit towards $-\infty$ and the upper limit towards $+\infty$ and this switching of the control register of the FPU considerably slows down the execution. A number of implementations bypass this problem by systematically subtracting $\epsilon$ from the lower bound and adding $\epsilon$ to the upper bound. As harmless as this extra $\epsilon$ may seem, it feeds the natural demon of interval arithmetic, which constantly increases the size of intervals. The result is then certainly fair, but useless. The recent agreement on a standard will somewhat facilitate the use of the intervals~\cite{7140721,8277144}. However, it is not enough to couple a library~\cite{6893849,6893854} to its program to automatically benefit from interval calculation. Worse, beyond the transformation of types, literal constants and, more importantly, input-output, these are the expressions and algorithms that must be radically revised to curb the expansion of the resulting intervals.
This technique has already been used for a thesis in the context of calculating uncertainties for data acquisitions~\cite{PhDFelixHautot}.
\subsection{Good news on accuracy front}\label{subsub:fma}
Since \textos{2008}~\cite{4610935} the IEEE~754 standard takes into account the operations \texttt{fma}~``\,\emph{Fused Multiply and Add}\,'' which allows one to perform only a rounding instead of two on the combination of a multiplication and an addition.

This new operation primarily affects two major algorithms: the scalar product~\cite{doi:10.1137/030601818} (and therefore the products of matrices but also convolutions) and the evaluation of polynomials by the \textsc{Horner-Ruffini} method.

These two algorithms, based on sums, are intended to be coupled with compensation techniques such as the summation algorithm of \textsc{Kahan}~\cite{CompensationLanglois,CompensationLouvet1,CompensationLouvet2}.

\def\rnd{\mathop{\mathgroup\symoperators{}rnd}\nolimits}%
\texttt{fma} can also be used for Error-Free Transformation (EFT) of arithmetics: given one operation $\circ$ and two floating point numbers $a$ and $b$, the rounded operation $a\circ{}b$ would produce an approximate result $s$, but the corresponding EFT will produce a sum of two floating points values $s+e$ corresponding to the exact result. This capacity is relevant to use when computing a reduced discriminant $\Delta'={b'^{2}-ac}$ riddled with cancellation when $b'^{2}\sim{}ac$. The EFT for the multiplication has two steps: the most significant result is computed with a simple, approximate (rounded) product $a*{}b$, also noted $a\otimes{b}$ distinct from the exact product $a\times{}b$, and the error term is evaluated with a \texttt{fma} as follows:
\begin{eqnarray}
a\times{}b = a\otimes{b} + \mathtt{fma}(a,b,-a\otimes{b})
\end{eqnarray}
One then apply this EFT to the two terms in $\Delta'={b'^{2}-ac}$ and if there is a cancellation between the most significant terms of the EFT, adding the difference of the least significant terms will provide for the missing accuracy.

\subsection{Stochastic Arithmetic}
Beyond these self-defense and stability validation approaches, we can use more recent environments for assessing digital precision. These execution environments induce small disturbances at each calculation stage and thus make it possible to measure the fluctuation of the result. Here the emphasis is about stochastic arithmetic. This technique is akin to Monte Carlo. The disturbance induced in the computation can appear in the form of either a random noise added to the computation or a random change of mode of rounding of the results.
These environments are \prog{CADNA}~\cite{SCOTT2007507,JEZEQUEL2008933,JEZEQUEL20101927,LAMOTTE20101925,CADNA}, \prog{Verificarlo}~\cite{Verificarlo,denis2015verificarlo} and \prog{Verrou}~\cite{fevotte2016verrou,fevotte:hal-01908466,fevotte:hal-01383417,Verrou} which can be used with \prog{valgrind}~\cite{refProgramValgrind,NETHERCOTE200344,10.1145/1273442.1250746}.

The \prog{Verrou} environment allows, in addition to its stochastic mode, to force all rounding of the calculation in an unconventional direction instead of the nearest rounding: towards $-\infty$, towards $+\infty$, towards $0$, and even to test the rounding off ``as far as possible'', a good way to get an idea of the resilience of ones numerical code.
These programs allow an instrumentation of the code which reveals numerical weaknesses at runtime. Their implementation has an impact on the performance of the code, and must be taken into account in the time budget for testing and optimization.

\subsection{Double precision as an insurance}
The previous examples have shed a light on numerical computation that is far from common expectation: floating-point computing is not straightforward at all, while expressions like ``number crunching'' evoke blind use of power.
In fact, while a mathematical formula looks like a monolithic, deep-frozen archive of theoretical truth, its practical evaluation often requires to be split between different domains and regimes, each relying on mathematically identical but numerically different expressions. This casuistry may feel uncomfortable for the beginner, but it is also irritating to the optimization expert: as pipelining and vectorization are the modern basis of computing performance, tests in the innermost loops of library are shunned.
In the same way, cache-misses induced by lookup table are also a pain for performance computing.

All fields of engineering usually involve ``\emph{factors of safety}'' to provide for unexpected situation. The fields where the materials are known the better can afford to minimize these factors of safety: concrete walls or elevator cables can be oversized with security factors of order twenty, but planes have to stay within a security factor of two: this can only be achieved through a deep and thoroughly documented knowledge of the material at hand. One might see the traditional recourse to double precision as such a factor of safety. While most computations will not be used in real life with more than single precision, intermediate result may require double precision to achieve single accuracy. For end user to rely only on single precision, and possibly benefit on its gain in speed, they need a much higher level of test and documentation of their algorithms and their library implementation.

\subsection{Mixed precision}\label{sub:mixprecision}
As mentioned in Subsection~\ref{sub:minimization}, the algebraic intensity of machine learning computations and the best use of both memory size and bandwidth had the community turn to accuracies lower than single precision.
Namely half-precision, either the IEEE compliant \texttt{binary16}, or \texttt{FP16}~\cite{4610935}, or the \texttt{bfloat16} of Google Brain~\cite{DBLP:journals/corr/abs-1905-12322}.
\texttt{FP16} having a very restricted range (magnitude up to $65\,536$), one often use a single precision accumulator while multiplying half precision variable.

Besides machine learning, dense and sparse linear algebra, as well as Fast Fourier Transform are studied with mixed precision in Ref.~\cite{abdelfattah2020survey}.

\subsection{Software engineering}\label{sub:software_engineering}
The major programming paradigms prevents us from dangerous practices through discipline:
structured programming forbids the use of \texttt{goto},
modular programming prevents global variables,
functional programming forbids mutations,
strong-typing prevents implicit cast of type,
polymorphism discipline the indirect transfer of control of pointers to functions,~\emph{etc}.
We should be looking for what to forbid to achieve a safer floating-point environment:
let us note that it has not to rely on a specific syntax or language, but more on good practices.

The first and most well known approach is that
tests for strict equalities should likely be banned (topological aspect of the representation),
but it is only the most elementary defense.


To coin Robert \textsc{Martin}~\cite{martin2018clean} 
``A good architecture allows you to defer framework decisions''. 
The choice of the precision is such a framework decision. In order to defer this choice, \emph{generic programming\ } is the paradigm of choice. The easiest way to implement it, is using simple preprocessing of source or \texttt{C} style \texttt{typedef} to specify at compilation time a generic floating-point type with a given precision.
Unfortunately, the precision of the floating-point variables does not guarantee in itself the accuracy of the result. Hence the full generic paradigm is required here, allowing to rewrite the function differently for the various precision if needed.


Unit testing and Test Driven Development may also look beneficial to the numerical computation community, but the amount of testing required is rather large.

\section{Conclusion}
    {``\emph{The purpose of computing is insight, not numbers.}'' ({Richard Wesley \textsc{Hamming}, Turing Award \textos{1968}}})

Experimental systematic error usually represent a large part of the discussion in a science paper or thesis. However numerical systematic error are rarely discussed.
The propagation of numerical errors along the computational process needs to be given attention as part of the systematic error study.

Floating-point numbers are a computer abstraction for real numbers. As we have shown in this paper, this abstraction does not protect us completely from implementation details: to coin Ref.~\cite{JohnDCook,LeakyAbstractions} floating-point numbers are a \emph{leaky abstraction}.

Fifteen years after the computer power growth is no more being driven by the processor frequency, in order to achieve the expectation on science projects, we face the need for keeping the computing power growing, and some in the field consider moving to less precise computation to execute it faster.

Admittedly double precision is not double accuracy. Floating-point precision is a mean to achieve floating-point accuracy, the goal of our computation, that has to be asserted beforehands. The problem stands in computations larger than the ones we can follow, where not only the accumulation of round-off errors, but the existence of numerical singularities hidden in the more and more complex domain spaces of these computation hinders our confidence in floating points computations.
In the same way the software engineering has dealt with the growing complexity of modern programs over seventy years, through a set of disciplines enforced by language paradigms as well as helper tools, floating point computing must be understood as the difficult field it is, deserving also its set of good practices and tools.

As W.~\textsc{Kahan} reminds us in the conclusion of Ref. \cite{kahan1996beastly}, we still have to rely on the usual strategies: different compilers, and if possible different computers (not only a different hardware, but a different architecture, and possibly a different manufacturer), different libraries and different runtime environments also matter.
Using binary format data storage from end to end is also a good approach, separating the data model from the data view.

The self-defense techniques suggest us to develop our codes in generic precision, instead of setting from start a fixed precision. The relative safety provided by double precision is dismissed by teams ready for more aggressive approach, particularly with large amount of data, naively confusing precision with accuracy. An empirical, down to earth approach, comparing results of different precisions, should help decide on the precision. Unfortunately, conducting such tests is not as simple as it seems~\cite{Monniaux:2008:PVF:1353445.1353446}.

\begin{acknowledgements}
We gratefully acknowledge the stimulating discussions we had with Amel Korichi, Hadrien Grasland, Pierre Aubert and Philippe Gauron.
\end{acknowledgements}

\bibliographystyle{spphys}
\Urlmuskip=0mu plus 2mu\relax
\bibliography{bibliographie}

\begin{thebibliography}{10}
\providecommand{\url}[1]{{#1}}
\providecommand{\urlprefix}{URL }
\expandafter\ifx\csname urlstyle\endcsname\relax
  \providecommand{\doi}[1]{DOI \discretionary{}{}{}#1}\else
  \providecommand{\doi}{DOI \discretionary{}{}{}\begingroup
  \urlstyle{rm}\Url}\fi

\bibitem{Goldberg:1991:CSK:103162.103163}
D.~Goldberg, ACM Comput. Surv. \textbf{23}(1), 5, Mar. {\oldstyle 1991}.
\newblock \doi{10.1145/103162.103163}

\bibitem{Kahan:1979:PFS:1057520.1057522}
W.~Kahan, J.F. Palmer, SIGNUM Newsl. \textbf{14}(si-2), 13, Oct. {\oldstyle
  1979}.
\newblock \doi{10.1145/1057520.1057522}

\bibitem{30711}
{IEEE Standard for Binary Floating-Point Arithmetic}.
\newblock Tech. rep., Oct {\oldstyle 1985}.
\newblock \doi{10.1109/IEEESTD.1985.82928}

\bibitem{27840}
{IEEE Standard for Radix-Independent Floating-Point Arithmetic}.
\newblock Tech. rep., Oct {\oldstyle 1987}.
\newblock \doi{10.1109/IEEESTD.1987.81037}

\bibitem{4610935}
{IEEE Standard for Floating-Point Arithmetic}.
\newblock Tech. rep., Aug {\oldstyle 2008}.
\newblock \doi{10.1109/IEEESTD.2008.4610935}

\bibitem{8766229}
{IEEE Standard for Floating-Point Arithmetic}.
\newblock Tech. rep., July {\oldstyle 2019}.
\newblock \doi{10.1109/IEEESTD.2019.8766229}

\bibitem{IEC559:1989}
{Binary floating-point arithmetic for microprocessor systems}.
\newblock Standard, International Organization for Standardization, Geneva, CH,
  Feb. {\oldstyle 1989}.
\newblock \urlprefix\url{https://www.iso.org/standard/19706.html}

\bibitem{ISO60559:2011}
{Information technology — Microprocessor Systems — Floating-Point
  arithmetic}.
\newblock Standard, International Organization for Standardization, Geneva, CH,
  Jun. {\oldstyle 2011}.
\newblock \urlprefix\url{https://www.iso.org/standard/57469.html}

\bibitem{ISO60559}
{Information technology — Microprocessor Systems — Floating-Point
  arithmetic}.
\newblock Standard, International Organization for Standardization, Geneva, CH,
  May {\oldstyle 2020}.
\newblock \urlprefix\url{https://www.iso.org/standard/80985.html}

\bibitem{10.1145/62138.62141}
B.~Liskov, in \emph{Addendum to the Proceedings on Object-Oriented Programming
  Systems, Languages and Applications (Addendum)} (Association for Computing
  Machinery, New York, NY, USA, 1987), OOPSLA ’87, p. 17–34.
\newblock \doi{10.1145/62138.62141}

\bibitem{10.1007/3-540-47910-4_8}
B.~Liskov, J.M. Wing, in \emph{ECOOP' 93 --- Object-Oriented Programming}, ed.
  by O.M. Nierstrasz (Springer Berlin Heidelberg, Berlin, Heidelberg, 1993),
  pp. 118--141

\bibitem{10.1145/197320.197383}
B.H. Liskov, J.M. Wing, ACM Trans. Program. Lang. Syst. \textbf{16}(6),
  1811–1841, Nov. {\oldstyle 1994}.
\newblock \doi{10.1145/197320.197383}.
\newblock \urlprefix\url{https://doi.org/10.1145/197320.197383}

\bibitem{castiel:inria-00000567}
A.~Castiel, V.~Lef{\`e}vre, P.~Zimmermann, {Interstices}  {\oldstyle 2004}.
\newblock \urlprefix\url{https://hal.inria.fr/inria-00000567}.
\newblock \url{http://interstices.info/display.jsp?id=c\_5936}. See also
  \url{http://www.vinc17.org/research/slides/epao2001-03.pdf}

\bibitem{10.1117/12.505591}
C.~Daramy, D.~Defour, F.~de~Dinechin, J.M. Muller, in \emph{Advanced Signal
  Processing Algorithms, Architectures, and Implementations XIII}, vol. 5205,
  ed. by F.T. Luk. International Society for Optics and Photonics (SPIE, 2003),
  vol. 5205, pp. 458 -- 464.
\newblock \doi{10.1117/12.505591}

\bibitem{QuickDirtyToDecimal}
R.~Regan.
\newblock {Quick and Dirty Floating-Point to Decimal Conversion}.
\newblock
  \href{https://www.exploringbinary.com/quick-and-dirty-floating-point-to-decimal-conversion}{\texttt{https://www.exploringbinary.com\hspace{0pt}/quick-and-dirty-\hspace{0pt}floating-point-to-decimal-conversion}},
  Nov. {\oldstyle 2010}

\bibitem{Gay90correctlyrounded}
D.M. Gay, {Correctly Rounded Binary-Decimal and Decimal-Binary Conversions}.
\newblock Tech. Rep. 90-10, Numerical Analysis Manuscript AT\&T Bell
  Laboratories, Murray Hill, New Jersey 07974, Nov. {\oldstyle 1990}.
\newblock \urlprefix\url{https://ampl.com/REFS/rounding.pdf}

\bibitem{Loitsch:2010:PFN:1809028.1806623}
F.~Loitsch, SIGPLAN Not. \textbf{45}(6), 233, Jun. {\oldstyle 2010}.
\newblock \doi{10.1145/1809028.1806623}

\bibitem{10.1145/363067.363112}
I.B. Goldberg, Commun. ACM \textbf{10}(2), 105–106, Feb. {\oldstyle 1967}.
\newblock \doi{10.1145/363067.363112}

\bibitem{Normale04whendouble}
S.~Boldo, G.~Melquiond, {When double rounding is odd}.
\newblock Tech. Rep. Research Report No 2004-48, École Normale Supérieure
  Lyon, Nov. {\oldstyle 2004}.
\newblock
  \urlprefix\url{http://www.ens-lyon.fr/LIP/Pub/Rapports/RR/RR2004/RR2004-48.pdf}

\bibitem{doi:10.1137/0914050}
N.J. Higham, SIAM Journal on Scientific Computing \textbf{14}(4), 783,
  {\oldstyle 1993}.
\newblock \doi{10.1137/0914050}

\bibitem{Kahan:1965:PRR:363707.363723}
W.~Kahan, Commun. ACM \textbf{8}(1), 40, Jan. {\oldstyle 1965}.
\newblock \doi{10.1145/363707.363723}.
\newblock \url{http://doi.acm.org/10.1145/363707.363723}

\bibitem{10.1145/980175.980177}
J.~Mcnamee, ACM SIGSAM Bulletin \textbf{38}, 1, 03 {\oldstyle 2004}.
\newblock \doi{10.1145/980175.980177}

\bibitem{doi:10.1137/1037130}
T.O. Espelid, SIAM Review \textbf{37}(4), 603,  {\oldstyle 1995}.
\newblock \doi{10.1137/1037130}

\bibitem{babuska1968numerical}
I.~Babuška, in \emph{IFIP Congress (1)}, vol.~68 (1968), vol.~68, pp. 11--23

\bibitem{doi:10.1137/030601818}
T.~Ogita, S.M. Rump, S.~Oishi, SIAM Journal on Scientific Computing
  \textbf{26}(6), 1955,  {\oldstyle 2005}.
\newblock \doi{10.1137/030601818}

\bibitem{doi:10.1137/07068816X}
S.M. Rump, T.~Ogita, S.~Oishi, SIAM Journal on Scientific Computing
  \textbf{31}(2), 1269,  {\oldstyle 2009}.
\newblock \doi{10.1137/07068816X}

\bibitem{10.2307/107508}
W.G. Horner, Philosophical Transactions of the Royal Society of London
  \textbf{109}, 308,  {\oldstyle 1819}.
\newblock \doi{10.2307/107508}

\bibitem{Knuth:1962:EPC:355580.369074}
D.E. Knuth, Commun. ACM \textbf{5}(12), 595, Dec. {\oldstyle 1962}.
\newblock \doi{10.1145/355580.369074}

\bibitem{Graillat:2006:ICH:1141277.1141585}
S.~Graillat, P.~Langlois, N.~Louvet, in \emph{Proceedings of the 2006 ACM
  Symposium on Applied Computing} (ACM, New York, NY, USA, 2006), SAC '06, pp.
  1323--1327.
\newblock \doi{10.1145/1141277.1141585}

\bibitem{doi:10.1080/00401706.1962.10490022}
B.P. Welford, Technometrics \textbf{4}(3), 419,  {\oldstyle 1962}.
\newblock \doi{10.1080/00401706.1962.10490022}.
\newblock
  \urlprefix\url{https://www.tandfonline.com/doi/abs/10.1080/00401706.1962.10490022}

\bibitem{Kahan_Triangle}
W.~Kahan,   {\oldstyle 2014}.
\newblock
  \urlprefix\url{https://people.eecs.berkeley.edu/~wkahan/Triangle.pdf}.
\newblock \url{https://people.eecs.berkeley.edu/~wkahan/Triangle.pdf}

\bibitem{Kahan_Tetrahedron}
W.~Kahan,   {\oldstyle 2012}.
\newblock
  \urlprefix\url{https://people.eecs.berkeley.edu/~wkahan/VtetLang.pdf}.
\newblock \url{https://people.eecs.berkeley.edu/~wkahan/VtetLang.pdf}

\bibitem{byckling1973particle}
E.~Byckling, K.~Kajantie,   {\oldstyle 1973}

\bibitem{RUMP1988109}
S.M. Rump, in \emph{Reliability in computing}, ed. by R.E. Moore (Elsevier,
  1988), pp. 109--126.
\newblock \doi{https://doi.org/10.1016/B978-0-12-505630-4.50012-2}.
\newblock
  \urlprefix\url{http://www.sciencedirect.com/science/article/pii/B9780125056304500122}

\bibitem{loh2002rump}
E.~Loh, G.W. Walster, Reliable Computing \textbf{8}(3), 245,  {\oldstyle 2002}

\bibitem{Kahan:CSD-92-667}
W.~Kahan, J.W. Thomas, Augmenting a programming language with complex
  arithmetic.
\newblock Tech. Rep. UCB/CSD-92-667, EECS Department, University of California,
  Berkeley, Dec {\oldstyle 1991}.
\newblock
  \urlprefix\url{http://www2.eecs.berkeley.edu/Pubs/TechRpts/1991/6127.html}

\bibitem{vanOldenborgh:1990yc}
G.J. van Oldenborgh, Comput. Phys. Commun. \textbf{66}, 1,  {\oldstyle 1991}.
\newblock \doi{10.1016/0010-4655(91)90002-3}

\bibitem{courbariaux2015}
M.~Courbariaux, Y.~Bengio, J.P. David, in \emph{Advances in neural information
  processing systems} (2015), pp. 3123--3131.
\newblock \urlprefix\url{https://arxiv.org/abs/1412.7024}.
\newblock Accepted as a workshop contribution at ICLR 2015

\bibitem{DBLP:journals/corr/abs-1710-03740}
P.~Micikevicius, S.~Narang, J.~Alben, G.F. Diamos, E.~Elsen, D.~Garc{\'{\i}}a,
  B.~Ginsburg, M.~Houston, O.~Kuchaiev, G.~Venkatesh, H.~Wu, CoRR
  \textbf{abs/1710.03740},  {\oldstyle 2017}.
\newblock \urlprefix\url{http://arxiv.org/abs/1710.03740}

\bibitem{gupta2015deep}
S.~Gupta, A.~Agrawal, K.~Gopalakrishnan, P.~Narayanan, in \emph{International
  Conference on Machine Learning} (2015), pp. 1737--1746

\bibitem{DBLP:journals/corr/abs-1809-00095}
Y.~Choi, M.~El{-}Khamy, J.~Lee, CoRR \textbf{abs/1809.00095},  {\oldstyle
  2018}.
\newblock \urlprefix\url{http://arxiv.org/abs/1809.00095}

\bibitem{10.1093/qjmam/1.1.287}
A.M. Turing, The Quarterly Journal of Mechanics and Applied Mathematics
  \textbf{1}(1), 287, 01 {\oldstyle 1948}.
\newblock \doi{10.1093/qjmam/1.1.287}

\bibitem{7140721}
{IEEE Standard for Interval Arithmetic}.
\newblock Tech. rep., June {\oldstyle 2015}.
\newblock \doi{10.1109/IEEESTD.2015.7140721}.
\newblock \url{https://ieeexplore.ieee.org/document/7140721}

\bibitem{8277144}
{IEEE Standard for Interval Arithmetic (Simplified)}.
\newblock Tech. rep., Jan {\oldstyle 2018}.
\newblock \doi{10.1109/IEEESTD.2018.8277144}

\bibitem{6893849}
\emph{{2014 IEEE Conference on Norbert Wiener in the 21st Century (21CW)}}
  (Wiley-IEEE Computer Society Press, 2014).
\newblock \doi{10.1109/NORBERT.2014.6893849}

\bibitem{6893854}
M.~{Nehmeier}, in \emph{2014 IEEE Conference on Norbert Wiener in the 21st
  Century (21CW)} \cite{6893849}, pp. 1--6.
\newblock \doi{10.1109/NORBERT.2014.6893854}.
\newblock \url{https://github.com/nehmeier/libieeep1788} Its main focus is the
  correctness and not the performance of the implementation.

\bibitem{PhDFelixHautot}
F.~Hautot, {Cartographie topographique et radiologique 3D en temps réel :
  acquisition, traitement, fusion des données et gestion des incertitudes}.
\newblock Ph.D. thesis, Université Paris Saclay,  {\oldstyle 2017}

\bibitem{CompensationLanglois}
S.~Graillat, P.~Langlois, N.~Louvet.
\newblock {FMA Implementations of the Compensated Horner Scheme}, Jan.
  {\oldstyle 2006}.
\newblock \urlprefix\url{https://perso.univ-perp.fr/langlois
  /slides/dag06_sl.pdf}.
\newblock Reliable Implementation of Real Number Algorithms: Theory and
  Practice, Dagstuhl Seminar

\bibitem{CompensationLouvet1}
P.~Langlois, N.~Louvet.
\newblock {Operator Dependant Compensated Algorithms}, May {\oldstyle 2007}.
\newblock \urlprefix\url{http://perso.ens-lyon.fr/nicolas.louvet/LaLo07b.pdf}

\bibitem{CompensationLouvet2}
S.~Graillat, P.~Langlois, N.~Louvet.
\newblock {Accurate dot products with FMA}, Jul. {\oldstyle 2006}.
\newblock \urlprefix\url{http://rnc7.loria.fr/louvet_poster.pdf}

\bibitem{SCOTT2007507}
N.~Scott, F.~Jézéquel, C.~Denis, J.M. Chesneaux, Computer Physics
  Communications \textbf{176}(8), 507,  {\oldstyle 2007}.
\newblock \doi{https://doi.org/10.1016/j.cpc.2007.01.005}

\bibitem{JEZEQUEL2008933}
F.~Jézéquel, J.M. Chesneaux, Computer Physics Communications
  \textbf{178}(12), 933,  {\oldstyle 2008}.
\newblock \doi{https://doi.org/10.1016/j.cpc.2008.02.003}

\bibitem{JEZEQUEL20101927}
F.~Jézéquel, J.M. Chesneaux, J.L. Lamotte, Computer Physics Communications
  \textbf{181}(11), 1927,  {\oldstyle 2010}.
\newblock \doi{https://doi.org/10.1016/j.cpc.2010.07.012}

\bibitem{LAMOTTE20101925}
J.L. Lamotte, J.M. Chesneaux, F.~Jézéquel, Computer Physics Communications
  \textbf{181}(11), 1925,  {\oldstyle 2010}.
\newblock \doi{https://doi.org/10.1016/j.cpc.2010.07.006}

\bibitem{CADNA}
N.~Scott, F.~Jézéquel, C.~Denis, J.M. Chesneaux.
\newblock {The CADNA library},  {\oldstyle 2007}.
\newblock \urlprefix\url{http://cadna.lip6.fr/}

\bibitem{Verificarlo}
C.~Denis, P.D.O. Castro, E.~Petit.
\newblock {A tool for automatic Montecarlo Arithmetic analysis.}
\newblock \urlprefix\url{https://github.com/verificarlo/verificarlo}

\bibitem{denis2015verificarlo}
C.~Denis, P.D.O. Castro, E.~Petit.
\newblock Verificarlo: checking floating point accuracy through monte carlo
  arithmetic, Sep. {\oldstyle 2015}

\bibitem{fevotte2016verrou}
F.~F{\'e}votte, B.~Lathuili{\`e}re, SCAN 2016 p.~47,  {\oldstyle 2016}

\bibitem{fevotte:hal-01908466}
F.~F{\'e}votte, B.~Lathuili{\`e}re, in \emph{{13e colloque national en calcul
  des structures}} ({Universit{\'e} Paris-Saclay}, Giens, Var, France, 2017).
\newblock \urlprefix\url{https://hal.archives-ouvertes.fr/hal-01908466}

\bibitem{fevotte:hal-01383417}
F.~F{\'e}votte, B.~Lathuili{\`e}re, {VERROU: Assessing Floating-Point Accuracy
  Without Recompiling}, Oct. {\oldstyle 2016}.
\newblock \urlprefix\url{https://hal.archives-ouvertes.fr/hal-01383417}.
\newblock Working paper or preprint

\bibitem{Verrou}
F.~F{\'e}votte, B.~Lathuili{\`e}re.
\newblock {Verrou: a floating-point rounding errors checker}.
\newblock \urlprefix\url{http://edf-hpc.github.io/verrou/vr-manual.html}.
\newblock Verrou is french for ``lock''

\bibitem{refProgramValgrind}
 \urlprefix\url{http://www.valgrind.org/}

\bibitem{NETHERCOTE200344}
N.~Nethercote, J.~Seward, Electronic Notes in Theoretical Computer Science
  \textbf{89}(2), 44 ,  {\oldstyle 2003}.
\newblock \doi{10.1016/S1571-0661(04)81042-9}.
\newblock RV '2003, Run-time Verification (Satellite Workshop of CAV '03)

\bibitem{10.1145/1273442.1250746}
N.~Nethercote, J.~Seward, SIGPLAN Not. \textbf{42}(6), 89–100, Jun.
  {\oldstyle 2007}.
\newblock \doi{10.1145/1273442.1250746}

\bibitem{DBLP:journals/corr/abs-1905-12322}
D.D. Kalamkar, D.~Mudigere, N.~Mellempudi, D.~Das, K.~Banerjee, S.~Avancha,
  D.T. Vooturi, N.~Jammalamadaka, J.~Huang, H.~Yuen, J.~Yang, J.~Park,
  A.~Heinecke, E.~Georganas, S.~Srinivasan, A.~Kundu, M.~Smelyanskiy, B.~Kaul,
  P.~Dubey, CoRR \textbf{abs/1905.12322},  {\oldstyle 2019}.
\newblock \urlprefix\url{http://arxiv.org/abs/1905.12322}

\bibitem{abdelfattah2020survey}
A.~Abdelfattah, H.~Anzt, E.G. Boman, E.~Carson, T.~Cojean, J.~Dongarra,
  M.~Gates, T.~Grützmacher, N.J. Higham, S.~Li, N.~Lindquist, Y.~Liu, J.~Loe,
  P.~Luszczek, P.~Nayak, S.~Pranesh, S.~Rajamanickam, T.~Ribizel, B.~Smith,
  K.~Swirydowicz, S.~Thomas, S.~Tomov, Y.M. Tsai, I.~Yamazaki, U.M. Yang,
  {\oldstyle 2020}

\bibitem{martin2018clean}
R.C. Martin, \emph{Clean Architecture: A Craftsman's Guide to Software
  Structure and Design}.
\newblock Martin, Robert C (Prentice Hall, 2018).
\newblock \urlprefix\url{https://books.google.fr/books?id=8ngAkAEACAAJ}

\bibitem{JohnDCook}
J.D. Cook.
\newblock {Floating point numbers are a leaky abstraction},  {\oldstyle 2009}.
\newblock
  \urlprefix\url{https://www.johndcook.com/blog/2009/04/06/numbers-are-a-leaky-abstraction/}

\bibitem{LeakyAbstractions}
J.~Spolsky.
\newblock {The Law of Leaky Abstractions},  {\oldstyle 2002}.
\newblock
  \urlprefix\url{https://www.joelonsoftware.com/2002/11/11/the-law-of-leaky-abstractions/}

\bibitem{kahan1996beastly}
W.~Kahan.
\newblock Beastly numbers, Jan. {\oldstyle 1996}.
\newblock
  \urlprefix\url{https://people.eecs.berkeley.edu/~wkahan/tests/numbeast.pdf}

\bibitem{Monniaux:2008:PVF:1353445.1353446}
D.~Monniaux, ACM Trans. Program. Lang. Syst. \textbf{30}(3), 12:1, May
  {\oldstyle 2008}.
\newblock \doi{10.1145/1353445.1353446}

\end{thebibliography}
\end{document}